\documentclass[a4paper]{amsart}

\usepackage{amsthm}
\usepackage{amsmath, amssymb}
\usepackage[all]{xy}
\SelectTips{cm}{}

\newtheorem{theorem}{Theorem}[section]
\newtheorem{lemma}[theorem]{Lemma}
\newtheorem{proposition}[theorem]{Proposition}
\newtheorem*{theorem*}{Theorem}
\newtheorem{corollary}[theorem]{Corollary}
\theoremstyle{remark}
\newtheorem{definition}[theorem]{Definition}

\newtheorem{remark}[theorem]{Remark}
\newtheorem{example}[theorem]{Example}
\numberwithin{equation}{section}
\newcommand{\C}{C^*\text{-algebra}}
\newcommand{\W}{\text{Wojciech}}

\newcommand{\im}{\operatorname{im }}
\newcommand{\coker}{\operatorname{coker }}

\newcommand{\Z}{\mathbb{Z}}
\newcommand{\T}{\mathbb{T}}
\newcommand{\N}{\mathbb{N}}

\newcommand{\K}{\mathcal{K}}

\newcommand{\Aut}{\operatorname{Aut}}
\newcommand{\Prim}{\operatorname{Prim}}

\newcommand{\supp}{\operatorname{supp}}
\newcommand{\Hi}{\mathcal{H}}
\begin{document}
\title[The $C^*$-algebras of graphs with
sinks]{\boldmath{Classification theorems for the
$C^*$-algebras\\ of graphs with sinks}}
\author{Iain Raeburn}
\address{School of Mathematical and Physical Sciences\\ University of Newcastle\\
Callaghan\\NSW 2308\\ Australia} 
\email{iain@math.newcastle.edu.au}
\author{Mark Tomforde}
\address{Department of Mathematics\\ University of Iowa \\ Iowa City\\ Iowa 52242\\ USA}
\email{tomforde@math.uiowa.edu}
\author{Dana P. Williams}
\address{Department of Mathematics\\ Dartmouth College\\
Hanover\\ NH 03755\\ USA}
\email{dana.williams@dartmouth.edu}
\thanks{This research was supported by the Australian Research
Council.}
\begin{abstract}
We consider graphs $E$ which have been obtained by adding one or
more sinks to a fixed directed graph $G$. We classify the
$C^*$-algebra of $E$ up to a very strong equivalence relation,
which insists, loosely speaking, that $C^*(G)$ is kept fixed.
The main invariants are vectors $W_E:G^0\to \N$ which describe
how the sinks are attached to $G$; more precisely, the invariants
are the classes of the $W_E$ in the cokernel of the map $A-I$, where
$A$ is the adjacency matrix of the graph $G$.
\end{abstract}
\maketitle
The Cuntz-Krieger algebras $\mathcal{O}_A$ are generated by families of
partial isometries satisfying relations determined by a finite matrix
$A$ with entries in $\{0,1\}$ and no zero rows \cite{CK}. One can view
$\mathcal{O}_A$ as the
$C^*$-algebra of the finite directed graph $E$ with vertex adjacency
matrix
$A$ \cite{W}; note that $E$ has no sinks because $A$ has no zero rows. In
recent years there has been a flurry of interest in analogues of
these algebras for infinite graphs and matrices (see
\cite{KPR} and
\cite{EL}, for example). 
It was shown in \cite{RS} that the graph algebras of \cite{KPR}
and the Exel-Laca algebras of \cite{EL} can be realised as
direct limits of
$C^*$-algebras of finite graphs with sinks. Since sinks were
specifically excluded in the original papers, it is now of some
interest to investigate the effect of sinks on the structure of the graph
algebra and its $K$-theory. The results of \cite{BPRS} and \cite{RS}
show that this effect can be substantial, depending on how the sink is
attached to the rest of the graph. Here we shall prove some classification
theorems for graphs with sinks which describe the effect of
adding sinks to a given graph.
Suppose $E$ is a row-finite graph with one sink $v$. The set
$\{v\}$ is hereditary, and therefore gives rise to an ideal
$I(v)$ in the $C^*$-algebra $C^*(E)$ of $E$. According to general
theory, the quotient $C^*(E)/I(v)$ can be identified with the  graph
algebra $C^*(G)$ of the graph $G$ obtained, loosely speaking, by
deleting $v$ and all edges which head only into $v$ (see
\cite[Theorem~4.1]{BPRS}).  We consider primarily
graphs $E$ with one sink for which this quotient
is a fixed row-finite graph $G$; we call such graphs \emph{1-sink
extensions} of $G$ (see Definition~\ref{extension}).  The results in
\cite{RS} suggest that the appropriate invariant should be the 
\emph{$\W$ vector} of the extension, which is the element 
$W_E$ of $\prod_{G^0} \Z$ whose $w$th
entry is the number of paths in $E^1 \setminus G^1$ from $w$
to the sink.  
 We
now state our main theorem as it applies to the finite
graphs which give simple Cuntz-Krieger algebras. We denote by $A_G$
the vertex matrix of a graph $G$, in which
$A_G(w_1,w_2)$ is the number of edges in
$G$ from $w_1$ to $w_2$. For any row-finite graph $G$, $A_G$ is
a well-defined map on the direct product $\prod_{G^0}\Z$ and $A_G^t$ is
well-defined on
$\bigoplus_{G^0}
\Z$. 
\begin{theorem*}
Suppose that  $E_1$
and $E_2$ are 1-sink extensions of a finite transitive graph $G$.
\begin{enumerate}
\item
If $W_{E_1} - W_{E_2} \in \im (A_G-I)$, then there exist a
1-sink  extension $F$ of $G$ and  embeddings $\phi_i : C^*(F)
\rightarrow C^*(E_i)$ onto full corners of $C^*(E_i)$ such that the
following diagram commutes
$$ \xymatrix{C^{*}(F) \ar[rr]^{\phi_i}
\ar[dr]_{\pi_{F}} && C^{*}(E_{i})
    \ar[dl]^{\pi_{E_i}}  \\
& C^{*}(G). & \\ }
$$
\item
If there exist $F$ and $\phi_i$ as above,
and if $\ker (A_G^t - I) = \{ 0 \}$, then $W_{E_1} -
W_{E_2} \in \im (A_G-I)$.
\end{enumerate}
\end{theorem*}
While the invariants we are dealing with are $K$-theoretic in
nature, and the proof of part (2) uses
$K$-theory, we give constructive proofs of part (1) and of the other main
theorems.  Thus  we can
actually find the graph $F$.  For example, if
$G$ is given by
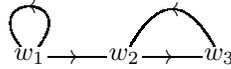
\begin{equation*}
\objectmargin{1pt}
  \xygraph{
!~:{@{-}|@{>}}
{w_{1}}(:@(ur,ul)"{w_{1}}"):[r] {w_{2}}:[r] {w_{3}} :@/_1.5pc/"{w_{2}}" }
\end{equation*}
and
$E_1$ and $E_2$ are the 1-sink extensions
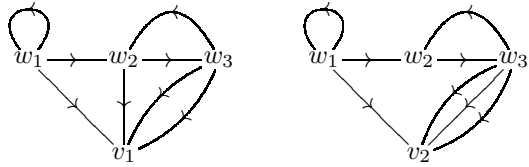
\begin{equation*}
\objectmargin{1pt}
  \xygraph{!~:{@{-}|@{>}}
{w_{1}}(:@(ur,ul)"{w_{1}}",:[dr]{v_{1}}):[r] {w_{2}} (:"{v_{1}}") :[r] {w_{3}}
   (:@/_1.5pc/"{w_{2}}", :@/^/"{v_{1}}", :@/_/"{v_{1}}") } \qquad
\xygraph{!~:{@{-}|@{>}}
{w_{1}}(:@(ur,ul)"{w_{1}}",:[dr]{v_{2}}):[r] {w_{2}} :[r] {w_{3}}
   (:@/_1.5pc/"{w_{2}}", :@/^.6pc/"{v_{2}}", :@/_.6pc/"{v_{2}}",
   :"{v_{2}}") } 
\end{equation*}
then we can take for $F$ the graph 
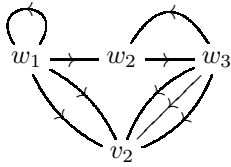
\begin{equation*}
  \xygraph{!~:{@{-}|@{>}}
{w_{1}}(:@(ur,ul)"{w_{1}}",:@/_/[dr]{v_{2}}, :@/^/"{v_{2}}")
  :[r] {w_{2}} :[r] {w_{3}}
   (:@/_1.5pc/"{w_{2}}", :@/^.6pc/"{v_{2}}", :@/_.6pc/"{v_{2}}",
   :"{v_{2}}") } 
\end{equation*}
The concrete nature of these constructions is very helpful when we want
to apply them to graphs with more than one sink, as we do in \S\ref{essn}.
It also means that our classification is quite different in
nature from the $K$-theoretic classifications of the algebras of
finite graphs without sinks \cite{Ro,H}. It
would be an interesting and possibly very hard problem to combine our
theorems with those of \cite{Ro, H} to say something about 1-sink
extensions of different graphs. 
\smallskip
We begin in \S\ref{sec-prelim} by establishing conventions and
notation. We give careful
definitions of 1-sink and
$n$-sink extensions, and describe the basic constructions which we
use throughout. In
\S\ref{ess1}, we consider a class of extensions which we call
\emph{essential}; these are the 1-sink extensions $E$ for which the ideal
$I(v)$ is an essential ideal in
$C^*(E)$. For essential 1-sink extensions of row-finite graphs we have a
very satisfactory classification (Theorem~\ref{essclass}), which
includes part (1) of the above theorem. We show by example that we cannot
completely discard the essentiality, but in 
\S\ref{iness1} we extend the analysis to cover non-essential
extensions $E_1$ and $E_2$ for which the 
primitive ideal spaces $\Prim C^*(E_1)$ and $\Prim C^*(E_2)$ are
appropriately homeomorphic. This extra generality is crucial in
\S\ref{essn}, where we use our earlier results to prove a classification
theorem for extensions with
$n$ sinks (Theorem~\ref{essnclass}). In our last section, we
investigate the necessity of our hypothesis on the Wojciech vectors.
In particular, part (2) of the above theorem follows from
Corollary~\ref{corconverse}. 
\section{Sink extensions and the basic constructions}
\label{sec-prelim}
  A directed graph $G = (G^0,G^1,r,s)$ consists of
a countable set $G^0$ of vertices, a countable set $G^1$ of
edges, and maps $r,s:G^1 \rightarrow G^0$ which identify the
range and source of an edge.   A vertex
$v \in G^0$ is a \emph{sink} if
$s^{-1}(v)=\emptyset$, or a \emph{source} if $r^{-1}(v) =
\emptyset$; $G$ is \emph{row-finite} if each
vertex emits at most finitely many edges.  All graphs in this paper are
row-finite and directed, and unless we say otherwise, $G$ will stand for a
generic row-finite graph. In general, our notation should be consistent with
that of
\cite{BPRS} and \cite{KPR}.
If $G$ is a row-finite graph, a \emph{Cuntz-Krieger
G-family} in a $C^*$-algebra consists of mutually orthogonal
projections $\{ p_v : v \in G^0 \}$ and
partial isometries $\{ s_e : e \in G^1 \}$ which satisfy the
\emph{Cuntz-Krieger relations} 
$$s_e^* s_e = p_{r(e)} \ \text{for} \  e \in G^1 \ \
\text{and} \
\ p_v =
\sum_{ \{ e : s(e) =v \} } s_e s_e^* \ \text{whenever } v \in G^0
\text{ is not a sink.}$$
We denote by $C^*(G)=C^*(s_e,p_v)$ the $C^*$-algebra of the
graph $G$, which is generated by a universal Cuntz-Krieger
$G$-family $\{s_e,p_v\}$ (see \cite[Theorem 1.2]{KPR}).
\begin{definition}\label{extension}
An \emph{$n$-sink extension} of $G$ is a row-finite graph $E$
which contains $G$ as a subgraph and satisfies:
\begin{enumerate}
\item
$H:=E^0\setminus G^0$ is finite, contains no sources,
and contains exactly $n$ sinks.
\item
There are no loops in $E$ whose vertices lie in
$H$.
\item
If $e \in E^1 \setminus G^1$, then $r(e) \in H$.
\item
If $w$ is a sink in $G$, then $w$ is a sink in $E$.
\end{enumerate}
When we say $(E, v_i)$ is an $n$-sink extension of $G$, we mean
that
$v_1,\cdots v_n$ are the $n$ sinks outside $G^0$. We consistently
write
$H$ for $E^0\setminus G^0$ and $S$ for the set of sinks
$\{v_1,\cdots,v_n\}$ lying in $H$.
If $w\in H$, then there are at most finitely many paths from $w$
to a given sink
$v_i$. If there is one sink $v_1$ and exactly one path from
every
$w\in H$ to
$v_1$, we call $(E,v_1)$ a \emph{1-sink
tree extension} of $G$. Equivalently, $(H,s^{-1}(H))$ is a tree.
\end{definition}
If we start with a graph $E$ with $n$ sinks, these ideas should
apply as follows.  Let $H$ be the saturation of the set $S$
of sinks in the sense of \cite{BPRS}, and take $G:=E\setminus
H:=(E^0\setminus H, E^1\setminus r^{-1}(H))$. Then $E$ satisfies
all the above properties with respect to $G$ except possibly (1);
if, however,
$E$ is finite and has no sources, this is automatic too. So the
situation of Definition~\ref{extension} is quite general.
Property (4)  ensures that the saturation of $S$ does not
extend into $G$; it also implies that an
$m$-sink extension of an $n$-sink extension of $G$ is an
$(m+n)$-sink extension of $G$, which is important for an
induction argument in \S\ref{essn}.
 
\begin{lemma}
Let $(E,v_i)$ be an $n$-sink extension of $G$.  Then
$H := E^0
\setminus G^0$ is a saturated hereditary subset of $E^0$.
Indeed, $H$ is the saturation $\overline{ S }$
of $ S:=\{v_1,\cdots ,v_n\}$.
\label{saturated}
\end{lemma}
\begin{proof}
Property (3) of
Definition~\ref{extension} implies that $H$ is hereditary,
and property (4) that
$H$ is saturated. Because $\overline{S }$ is the smallest
saturated set containing $S$, it now suffices to prove that
$H \subset \overline{S}$.  Suppose
that $w \notin
\overline{ S}$.  Then either there is a
path $\gamma$ from
$w$ to a sink $r(\gamma)\notin S$, or there is
an  infinite path which begins at $w$. In the first case, $w$
cannot be in $H$ because $r(\gamma)\notin H$ and $H$ is
hereditary. In the second case, $w$ cannot be in $H$ because
otherwise we would have an infinite path going round the finite
set $H$, and there would have to be a loop in $H$. Either way,
therefore, $w\notin H$, and we have proved $H \subset \overline{
S }$.
\end{proof}
\begin{corollary}\label{defpi}
Suppose that $(E,v_i)$ is an $n$-sink extension of $G$, and
$I(S)$ is the ideal in $C^*(E)=C^*(s_e,p_v)$ generated by the
projections $p_{v_i}$ associated to the sinks $v_i\in S$. Then
there is a surjection
$\pi_E$ of
$C^*(E)$ onto $C^*(G)=C^*(t_f,q_w)$ such that $\pi_E(s_e)=t_e$
for $e\in G^1$ and $\pi_E(p_v)=q_v$ for $v\in G^0$, and $\ker
\pi_E=I(S)$.
\end{corollary}
\begin{proof}
From  Lemma~\ref{saturated} and \cite[Lemma~4.3]{BPRS}, we see
that $I(S)=I(H)$, and the result follows from
\cite[Theorem~4.1]{BPRS}.
\end{proof}
\begin{definition}  An $n$-sink extension $(E,v_i)$ of
$G$ is \emph{simple} if $E^0\setminus G^0 =
\{v_i,\cdots, v_n\}$.
We want to associate to each $n$-sink extension $(E,v_i)$ a
simple extension by collapsing paths which end at one of the
$v_i$. For the precise definition, we need some notation.  An
edge $e$ with
$r(e)\in H$ and $s(e)\in G^0$ is called a \emph{boundary edge};
the sources of these edges are called \emph{boundary vertices}.
We write $B_E^1$ and $B_E^0$ for the sets of boundary edges and
vertices. If $w\in G^0$ and $1\leq i\leq n$,
we denote by $Z(w,v_i)$ the set of paths $\alpha$ from $w$ to
$v_i$ which leave
$G$ immediately in the sense that $r(\alpha_1)\in H$. The
\emph{Wojciech vector} of the sink $v_i$ is the element
$W_{(E;v_i)}$ of
$\prod_{G^0}\N$ given by
\[
W_{(E;v_i)}(w):=\#Z(w,v_i)\ \text{for $w\in G^0$;}
\]
notice that $W_{(E;v_i)}(w)=0$ unless $w$ is a boundary vertex.
If $E$ has just one sink, we denote its only Wojciech vector by
$W_E$.
The \emph{simplification} of 
$(E,v_i)$ is the graph
$SE$ with $(SE)^0:=G^0 \cup \{v_1,\cdots,v_n \}$,
\begin{gather*}
(SE)^1 := G^1 \cup \{ e^{(w,\alpha)} : w \in B_E^0
\mbox{ and } \alpha
\in Z(w,v_i) \mbox{ for some $i$}\},\\
s|_{G^1} = s_E, \quad s(e^{(w,\alpha)})
= w, \quad r|_{G^1} = r_E, \quad\mbox{
and}\quad r(e^{(w,\alpha)}) = r(\alpha).
\end{gather*}
\end{definition}
The simplification of $(E,v_i)$ is a simple $n$-sink extension of
$G$ with the same Wojciech vectors as $E$. We now describe how
the graph algebras are related:
\begin{proposition}\label{simpembeds}  Let $(E,v_i)$ be an
$n$-sink extension of
$G$, and let $\{s_e, p_v \}$, $\{t_f, q_w \}$ denote the
canonical Cuntz-Krieger families in $C^*(SE)$
and $C^*(E)$.  Then there is an embedding $\phi^{SE}$ of
$C^*(SE)$ onto the full corner in
$C^*(E)$ determined by the projection $\sum_{i=1}^n q_{v_i} +
\sum\{q_w:w \in G^0\}$, which satisfies $\phi^{SE}(p_{v}) =
q_{v}$ for all $v\in G^0\cup\{v_i\}$, and for which we have
a commutative diagram involving the maps $\pi_E$ of
Corollary~\ref{defpi}:
\begin{equation*}
  \xymatrix{C^{*}(SE) \ar[rr]^{\phi^{SE}}
\ar[dr]_{\pi_{SE}} && C^{*}(E)
    \ar[dl]^{\pi_{E}}  \\
& C^{*}(G)}.
\end{equation*}
\label{reduction}
\end{proposition}
\begin{proof}
The elements
\[
P_v := q_v \quad \text{and} \quad S_e :=
\begin{cases} t_e & \text{if $e \in G^1$} \\ t_\alpha &
\text{if $e = e^{(w,\alpha)}$} \end{cases}
\]
form a Cuntz-Krieger
$(SE)$-family in $C^*(E)$, so there is a homomorphism
$\phi^{SE}:=\pi_{S,P} : C^*(SE) \rightarrow C^*(E)$ with
$\phi^{SE}(p_v)=P_v$ and
$\phi^{SE}(s_e) = S_e$. We trivially have
$\phi^{SE}(p_{v}) = q_{v}$ for $v\in G^0\cup S$.
To see that $\phi^{SE}$ is injective, we use the
universal property of $C^*(E)$ to build an action $\beta:\T\to
\Aut C^*(E)$ such that
\[\beta_z(q_w) = q_w \quad \text{and} \quad \beta_z(t_f) =
\begin{cases} zt_f & \text{if $s(f) \in G^0$} \\ t_f &
\text{otherwise,} \end{cases}
\]
note that $\phi^{SE}$ converts the gauge action on $C^*(SE)$ to
$\beta$, and apply the gauge-invariant uniqueness theorem
\cite[Theorem 2.1]{BPRS}.
 It follows from \cite[Lemma
1.1]{BPRS} that $\sum_{i=1}^n q_{v_i} +
\sum\{q_w:w \in G^0\}$ converges strictly to a projection
$q \in M(C^*(E))$ such that
\[q t_\alpha
t_\beta^* = \begin{cases} t_\alpha t_\beta^* &\text{if
$s(\alpha) \in G^0 \cup S$} \\ 0 &\text{otherwise}.\\
\end{cases}
\]
Thus $qC^*(E)q$ is spanned by the elements
$t_{\alpha}t_{\beta}^*$ with $s(\alpha)=s(\beta) \in G^0 \cup
S$, and by applying the Cuntz-Krieger relations we may
assume
$r(\alpha)=r(\beta) \in G^0 \cup S$ also, so that the
range of $\phi$ is the corner $qC^*(E)q$. To see that this
corner is full, suppose $I$ is an ideal containing $qC^*(E)q$. 
Then \cite[Lemma 4.2]{BPRS} implies that
$K:= \{ v : q_v \in I \}$ is a saturated hereditary subset
of $E^0$; since $K$ certainly contains
$ G^0
\cup S$, we deduce that
$K=E^0$.  But then $I=C^*(E)$ by \cite[Theorem
4.1]{BPRS}.  Finally, to see that the diagram commutes,
we just need to check that
$\pi_{SE}$ and $\pi_E \circ \phi^{SE}$ agree on generators.
\end{proof}
It is convenient to have a name for the situation described in
this proposition:
\begin{definition}
Suppose $(E,v_i)$ and $(F,w_i)$ are $n$-sink extensions of $G$.
We say that $C^*(F)$ is \emph{$C^*(G)$-embeddable in $C^*(E)$}
if there is an isomorphism $\phi$ of $C^*(F)=C^*(s_e,p_v)$ onto a
full corner in
$C^*(E)=C^*(t_f,q_w)$ such that $\phi(p_{w_i})=q_{v_i}$ for
all~$i$ and
$\pi_E\circ\phi=\pi_F:C^*(F)\to C^*(G)$. If $\phi$ is an
isomorphism onto $C^*(E)$, we say that $C^*(F)$ is
\emph{$C^*(G)$-isomorphic to
$C^*(E)$}. 
\end{definition}
Notice that if $C^*(F)$ is $C^*(G)$-embeddable in $C^*(E)$, then
$C^*(F)$ is Morita equivalent to $C^*(E)$ in a way which
respects the common quotient $C^*(G)$.
We now describe the basic construction by which we manipulate
the Wojciech vectors of graphs.
\begin{definition}Let $(E,v_i)$ be an $n$-sink extension of $G$,
and let $e$ be a boundary edge such that $s(e)$ is not a
source of $G$.  The \emph{outsplitting} of
$E$ by $e$ is the graph $E(e)$ defined by 
\begin{gather*}
E(e)^0 := E^0 \cup \{
v'
\};
\quad E(e)^1 := (E^1 \setminus \{ e \} ) \cup \{ e' \}
\cup \{f' : f \in E^1
\text{ and } r(f) =s(e) \}\\
(r,s)|_{E^1 \setminus \{ e \}} := (r_E,s_E); \quad
r(e'):=r_E(e),\ s(e'):=v'; \quad r(f'):=v',\ s(f'):=s_E(f).
\end{gather*}
In general, we call $E(e)$ a \emph{boundary outsplitting} of $E$.
\end{definition}
The following example might help fix the ideas:
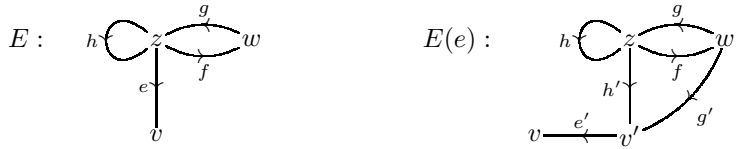
\begin{equation*}
\objectmargin{1pt}
E:\quad\xygraph{ 
 !~:{@{-}|@{>}}
 z (:@(ul,dl)"z"_{h},:[d]v_{e})  :@/_/[r] w_{f} :@/_/"z"_{g}
 }
\qquad\qquad\qquad E(e):\quad\xygraph{ 
 !~:{@{-}|@{>}}
 z (:@(ul,dl)"z"_{h},:[d]{v'}_{h'} :[l]v_{e'})  :@/_/[r] w_{f}
 (:@/_/"z"_{g}, :@/^/"{v'}"^{g'} ) }
\end{equation*}
If $(E,v_i)$ is an $n$-sink extension of
$G$,
then every boundary outsplitting $(E(e),v_i)$ is also an $n$-sink
extension of
$G$; if $(E,v_0)$ is a 1-sink tree extension, so is
$(E(e),v_0)$.  We need to assume that $s(e)$ is not a source of
$G$ to ensure that $E(e)$ is an $n$-sink extension, and we
make this assumption implicitly whenever we talk about boundary
outsplittings. As the name suggests, boundary outsplittings are
special cases of the outsplittings discussed in
\cite[\S2.4]{LM}.
\begin{proposition}\label{outspliteffectsW}
Suppose $(E(e),v_i)$ is a boundary outsplitting of an $n$-sink
extension
$(E,v_i)$  of
$G$. Then $C^*(E(e))$ is $C^*(G)$-isomorphic to $C^*(E)$. If $E$
is a 1-sink tree extension, then the Wojciech vector of $E(e)$ is
given in terms of the vertex matrix
$A_G$ of
$G$  by
\begin{equation}\label{outsplitw}
W_{E(e)} = W_E + (A_G-I) \ \delta_{s(e)}.
\end{equation}
\end{proposition}
\begin{proof} Let $C^*(E)=C^*(t_h, q_w)$. Then
\begin{align*}
P_v &:= \begin{cases} q_v & \text{if $v \neq s(e)$
and
$v
\neq v'$} \\ 
t_e t_e^* & \text{if $v =v'$} \\
q_{s(e)}-t_e t_e^* & \text{if $v = s(e)$} \\
\end{cases} \\S_g &:= \begin{cases} t_e & \text{if $g=e'$} \\ 
t_g(q_{s(e)}-t_e t_e^*) & \text{if $g
\neq e'$ and $r(g) = s(e)$}\\
t_f t_e t_e^* & \text{if $g = f'$ for some $f \in E^1$} \text{
with $r(f) = s(e)$}\\
t_g & \text{otherwise}
\end{cases} 
\end{align*}
is a Cuntz-Krieger $E(e)$-family which generates $C^*(E)$.  The
universal property of $C^*(E(e))=C^*(s_g,p_v)$ gives a
homomorphism
$\phi =\pi_{S,P}: C^*(E(e)) \rightarrow C^*(E)$ such that
$\phi (s_g) = S_g$ and $\phi (p_v) = P_v$, which is an
isomorphism by the gauge-invariant uniqueness theorem 
\cite[Theorem 2.1]{BPRS}. It is easy to check on generators that
$\phi$ is a $C^*(G)$-isomorphism.
When $H$ is a tree with one sink $v_{0}$, there is
precisely one path $\gamma$ in $E$ from $r(e)$ to $v_0$, and
hence all the new paths from a vertex $v$ to $v_0$ have the form
$f'\gamma$. Thus if
$v
\neq s(e)$,
\begin{align*}
W_{E(e)}(v)  
&=W_E(v) + \# \{ f' \in E(e)^1 : s(f')=v \text{ and }
f' \notin E^1 \} \\
&=W_E(v) + \# \{ f \in G^0 : s(f)=v \text{ and }
r(f)=s(e) \}\\
&=W_E(v) + A_G(v,s(e)).
\end{align*}
On the other hand, if $v=s(e)$, then
\begin{align*}
W_{E(e)}(s(e))  
&=W_E(s(e)) + \# \{ f' \in E(e)^1 : s(f')=s(e) \text{ and }
f' \notin E^1 \} - 1 \\
&=W_E(s(e)) + \# \{ f \in G^0 : s(f)=s(e)=
r(f) \} - 1\\
&=W_E(v) + A_G(s(e),s(e)) -1.
\end{align*}
Together these calculations give (\ref{outsplitw}).
\end{proof}
Suppose that $\alpha=\alpha_1\alpha_2\cdots\alpha_n$ is a path in
$G$ and there is a boundary edge
$e$ with $s(e)=r(\alpha)$.  Then  $E(e)$
will have a boundary edge
$\alpha_n'$ at $r(\alpha_{n-1})$, and therefore we can outsplit
again to get
$E(e)(\alpha_n')$. This graph has a boundary edge
$\alpha_{n-1}'$ at
$r(\alpha_{n-2})$, and we can outsplit again. Continuing this
process gives an extension
$E(e,\alpha)$ in which $s(\alpha)$ is a boundary
vertex.  We shall refer to this process as \emph{performing
outsplittings along the path $\alpha$.} From
Proposition~\ref{outspliteffectsW} we can calculate the Wojciech
vector of $E(e,\alpha)$:
\begin{corollary}
\label{outsplitalongpath}
Suppose $E$ is a 1-sink tree extension of $G$ and $\alpha$ is a
path in $G$ for which $r(\alpha)$ is a boundary vertex. Then for
any boundary edge $e$ with $s(e)=r(\alpha)$, we have
\[
W_{E(e,\alpha)}= W_E + \sum_{i=1}^{|\alpha|}
(A_G-I) \delta_{r(\alpha_i)}.
\]
\end{corollary}
\section{A classification for essential 1-sink
extensions}\label{ess1}
We now ask to what extent the $\W$ vector determines a 1-sink
extension. Suppose
that
$E_1$ and
$E_2$ are 1-sink extensions of $G$. Our main results say,
loosely speaking, that if the Wojciech vectors $W_{E_i}$
determine the same class in $\coker(A_G-I)$, then there will be
a simple extension $F$ such that $C^*(F)$ is
$C^*(G)$-embeddable in both 
$C^*(E_1)$ and $C^*(E_2)$. However, we shall need some
hypotheses on the way the sinks are attached to $G$; the
hypotheses in this section are satisfied if, for example, $G$ is
one of the finite transitive graphs for which $C^*(G)$ is a
simple Cuntz-Krieger algebra. We begin by describing the class of
extensions which we consider in this section. 
Recall that if
$v,w$ are vertices in $G$, then $v\geq w$ means there is a
finite path $\gamma$ with $s(\gamma)=v$ and $r(\gamma)=w$. For
$K, L \subset G^0$, $K \geq L$ means that for
each $v \in K$ there exists $w
\in L$ such that $v \geq w$.  If $\gamma$ is a loop, we
write $\gamma \geq L$ when $\{ r(\gamma_i) \}
\geq L$.  
\begin{definition}  A 1-sink extension $(E,v_0)$
of a graph $G$ is an \emph{essential extension} if $G^0\geq
v_0$. 
\end{definition}
We can see immediately that simplifications  of essential
extensions are essential, and consideration of a few cases shows
that boundary
outsplittings of essential extensions are essential. To see why we
chose the name, recall that an ideal 
$I$ in a $\C$ $A$ is  essential if $I \cap J \neq 0$
for all nonzero ideals $J$ in $A$, or equivalently, if $aI=0$
implies
$a=0$. Then we have:
\begin{lemma}  Let $(E, v_0)$ be a 1-sink extension of $G$. 
Then $(E,v_0)$ is an essential extension of $G$ if and only if
the ideal  
$I(v_0)$ generated by $p_{v_0}$ is an essential ideal in
$C^*(E)=C^*(s_e,p_v)$.
\label{border}
\end{lemma}
\begin{proof}
Suppose that
 there exists $w \in G^0$ such that
$w \ngeq v_0$. Then since
\[
I(v_0) = \overline{\text{span}} \{ s_\alpha s_\beta^* :
\alpha,\beta \in E^* \text{ and } r(\alpha)=r(\beta) = v_0 \}.
\]
(see \cite[Lemma 4.3]{BPRS}), we have $p_wI(v_0)=0$, and
$I(v_0)$ is not essential.  
Conversely, suppose that $G^0 \geq  v_0 $.    To show that
$I(v_0)$ is an essential ideal it suffices to prove that
if  $\pi : C^*(E) \rightarrow B(\Hi)$ is a
representation with $\ker
\pi \cap I(v_0) = \{ 0 \}$, then $\pi$ is faithful. 
So suppose  $\ker
\pi \cap I(v_0) = \{0 \}$. In particular, we have
$\pi(p_{v_0}) \neq 0$. For every $v\in G^0$  there is a
path $\alpha$ in $E$ such that $s(\alpha)=v$ and
$r(\alpha)=v_0$. Then
$\pi(s_\alpha^*s_\alpha) = \pi(p_{v_0}) \neq 0$, and hence
$\pi(p_v) \geq \pi(s_\alpha s_\alpha^*) \neq 0$.  Since every 
loop in a
1-sink extension $E$ must lie entirely in
$G$, every loop in $G$ has an exit in $E$; thus we can apply 
\cite[Theorem 3.1]{BPRS} to deduce that $\pi$ is faithful, as
required.
\end{proof}
We can now state our classification theorem for essential
extensions.
\begin{theorem}\label{essclass}
Let $G$ be a row-finite graph with no sources, and suppose that
$(E_1,v_1)$ and $(E_2,v_2)$ are essential 1-sink extensions of
$G$ with finitely many boundary vertices. If there exists
$n\in\bigoplus_{G^0}\Z$ such that the Wojciech vectors satisfy 
$W_{E_1}-W_{E_2}=(A_G-I)n$, then there is a simple 1-sink
extension $F$ of $G$ such that $C^*(F)$ is $C^*(G)$-embeddable
in both $C^*(E_1)$ and $C^*(E_2)$.
\end{theorem}
We begin by observing that, since a full corner in a full corner
of a $C^*$-algebra $A$ is a full corner in $A$, the composition
of two $C^*(G)$-embeddings is another $C^*(G)$-embedding. Thus it
suffices by  Proposition~\ref{simpembeds} to prove the theorem
for the simplifications
$SE_1$ and
$SE_2$. However, since we are going to perform boundary
outsplittings and these do not preserve simplicity, we assume
merely that $E_1$ and $E_2$ are 1-sink tree extensions. 
The following lemma is the key to many of our constructions:
\begin{lemma}
Let $(E_1,v_1)$ and $(E_2,v_2)$ be 1-sink tree extensions of
$G$ with finitely many boundary vertices, and suppose that
$B_{E_1}^0 \geq B_{E_2}^0\geq B_{E_1}^0$. If
$\gamma$ is a loop in $G$ such that $\gamma\geq B^0_{E_1}$, then
for any
$a
\in
\Z$ there are 1-sink tree extensions $E_1'$ and
$E_2'$ which are formed by performing a finite number of
boundary outsplittings to $E_1$ and $E_2$,
respectively, and for which
$$W_{E_1'} - W_{E_2'} = W_{E_1} - W_{E_2} +
a \Big(\sum_{j =1}^{|\gamma|} (A_G-I)
\delta_{r(\gamma_j)}\Big).$$
\label{Iainsideaforalemma}
\end{lemma}
\begin{proof}
Since the statement is symmetric in $E_1$ and $E_2$, it suffices to prove
this for $a>0$. Choose a path $\alpha$ in $G$ such that $s(\alpha)=r(\gamma)$
and $r(\alpha)\in B^0_{E_1}$. Since $B_{E_1}^0$ is finite, going
along paths from $r(\alpha)$ to $B_{E_2}^0$ and then to and fro
between $B_{E_2}^0$ and $B_{E_1}^0$ must eventually give either (a) a
loop $\mu$ which visits both $B_{E_1}^0$ and $B_{E_2}^0$, and a
path $\beta$ with $s(\beta)=r(\alpha)$ and $r(\beta)=s(\mu)\in
B_{E_1}^0$, or (b) a vertex $v\in B_{E_1}^0\cap B_{E_2}^0$ and a path $\beta$ with $s(\beta)=r(\alpha)$ and $r(\beta)=v$. 

We deal with case (a) first. Since there are boundary edges $e_1\in B_{E_1}^1$ and
$e_2\in B_{E_2}^1$ with $s(e_i)$ on $\mu$, we can perform
outsplittings along $\mu$ to get new tree extensions
$E_i(e_i,\mu^i)$, where
$\mu^i$ is the loop $\mu$ relabelled so that it ends at $s(e_i)$.
Because $\mu^1$ and $\mu^2$ have the same vertices as $\mu$ in a
different order, Corollary~\ref{outsplitalongpath} gives
\[
W_{E_i(e_i,\mu^i)}=W_{E_i} + \sum_{j=1}^{|\mu |} (A_G-I)
\delta_{r(\mu_j)},
\]
so we have $ W_{E_1(e_1,\mu^1)} -
W_{E_2(e_2,\mu^2)} = W_{E_1}-W_{E_2}$. Since
$r(\beta_{|\beta|})=s(\mu)$, and in forming both $E_i(e_i,\mu^i)$
we have performed an outsplitting at $s(\mu)$,
$s(\beta_{|\beta|})$ is a boundary vertex in both
$E_i(e_i,\mu^i)$; say $f_i\in B_{E_i}^1$ has
$s(f_i)=s(\beta_{|\beta|})$. Write
$\beta=\beta'\beta_{|\beta|}$, $\gamma^a$ for the path obtained by going
$a$ times around $\gamma$, and define
\[
E_1':=E_1(e_1,\mu^1)(f_1,\gamma^a\alpha\beta')\ \mbox{ and }\  
E_2':=E_2(e_2,\mu^2)(f_2,\alpha\beta').
\]
We now compute the Wojciech vectors using
Corollary~\ref{outsplitalongpath}: for example,
\[
W_{E_1'}=W_{E_1(e_1,\mu^1)}+(A_G-I)
\Big(\sum_{j=1}^{|\beta|-1}\delta_{r(\beta_j)}
+\sum_{j=1}^{|\alpha|}\delta_{r(\alpha_j)}
+\sum_{j=1}^{|\gamma|}a\delta_{r(\gamma_j)}\Big).
\]
The formula for $W_{E_2'}$ is the same except for the last
term, so
\begin{align*}
W_{E_1'}-W_{E_2'}&=W_{E_1(e_1,\mu^1)} -
W_{E_2(e_2,\mu^2)}
+\sum_{j=1}^{|\gamma|}a(A_G-I)\delta_{r(\gamma_j)}\\
&=W_{E_1}-W_{E_2}
+\sum_{j=1}^{|\gamma|}a(A_G-I)\delta_{r(\gamma_j)},
\end{align*}
as required.

In case (b), we can dispense with the first step in the preceding argument: we choose boundary edges $f_i\in B_{E_i}^1$ with $s(f_i)=v$, and then
\[
E_1':=E_1(f_1,\gamma^a\alpha\beta) \ \mbox{ and }\  
E_2':=E_2(f_2,\alpha\beta)
\]
have the required properties.
\end{proof}
\begin{proof}[Proof of Theorem~\ref{essclass}] 
As we indicated earlier, it suffices to prove the theorem when
$E_1$ and $E_2$ are tree extensions. It also suffices to prove
that we can perform boundary outsplittings on
$E_1$ and
$E_2$ to achieve extensions
$F_1$ and
$F_2$ with the same Wojciech vector;
Propositions~\ref{simpembeds} and
\ref{outspliteffectsW} then imply that we can take for $F$ the
common  simplification of $F_1$ and $F_2$. We can write
$n=\sum_{k=1}^m a_k \delta_{w_k}$ for some finite set $\{
w_1,w_2, \ldots , w_{m} \} \subset G^0$.  We shall prove by
induction on $m$ that we can perform the required outsplittings.
If $m=0$, then $W_{E_1} = W_{E_2}$, and there is
nothing to prove. So we 
suppose that we can perform the outsplittings
whenever $n$ has the form $\sum_{k=1}^m a_k \delta_{w_k}$,
and that
$n=\sum_{k=1}^{m+1} a_k \delta_{w_k}$. Let $D$ be the subgraph
of $G$ with vertices $D^0 := \{ w_1,w_2, \ldots , w_{m+1} \}$
and edges $D^1 := \{ e \in G^1 : s(e),r(e) \in D^0
\}$.  Since $D$ is a finite graph it contains either a sink or a
loop.
If $D$ contains a sink, then by relabelling we can assume the
sink is 
$w_{m+1}$. Since  $A_G(w_{m+1},w_j)=0$ for all $j$, we have 
\[
W_{E_1}(w_{m+1}) = 
W_{E_2}(w_{m+1}) - a_{m+1}. 
\]
Thus either $E_1$ or $E_2$ has at least
$|a_{m+1}|$ boundary edges leaving $w_{m+1}$: we may as well
assume that
$a_{m+1} > 0$, so that
$W_{E_2}(w_{m+1}) \geq a_{m+1}$.  We can then perform
$a_{m+1}$ boundary outsplittings on $E_2$ at
$w_{m+1}$ to get a new extension $E_2'$.  From
Proposition~\ref{outspliteffectsW}, we have $W_{E_2'} =
W_{E_2} + a_{m+1}(A_G-I) \delta_{w_{m+1}}$, and therefore
\[
W_{E_1} =W_{E_2'} + (A_G-I) \Big( \sum_{k=1}^m a_k
\delta_{w_k} \Big).
\]
Since $E_2'$ is formed by performing boundary outsplittings to
the essential tree extension $E_2$, it  is also an essential tree
extension, and the inductive hypothesis implies that
we can perform boundary outsplittings on $E_1$ and
$E_2'$ to arrive at extensions with  the same Wojciech vector.
If $D$ does not have a sink,  it must contain a loop
$\gamma$.  If necessary, we can
shrink $\gamma$ so that its vertices are distinct, and by
relabelling, we may assume that $w_{m+1}$ lies on $\gamma$.
Because the extensions are essential, we have $G^0\geq
B_{E_1}^0$ and
$G^0\geq B_{E_2}^0$, so we can apply
Lemma~\ref{Iainsideaforalemma}. Thus
there are 1-sink tree extensions $E_1'$ and $E_2'$ formed by
performing boundary outsplittings to $E_1$ and $E_2$, and for
which
\[
W_{E_1'}-W_{E_2'} = W_{E_1}-W_{E_2} - a_{m+1}
\sum_{j=1}^{|\gamma|} (A_G-I) \delta_{r(\gamma_j)}.
\]
But because $W_{E_1} = W_{E_2} + (A_G-I)n$ this implies
that 
\[
W_{E_1'}=W_{E_2'} + (A_G-I) \Big( \sum_{j=1}^{m} b_j
\delta_{w_j} \Big),
\]
where $b_j=a_j-a_{m+1}$ if $w_j$ lies on $\gamma$, and
$b_j=a_j$ otherwise.
We can now invoke the
inductive hypothesis to see that we can perform boundary
outsplittings to
$E_1'$ and
$E_2'$  to arrive at extensions with the
same Wojciech vector.
This completes the proof of the inductive step, and the result
follows.
\end{proof}
\begin{remark}\label{standcons}
The graph $F$ in Theorem~\ref{essclass} has actually been
constructed in a very specific way, and it will be important in
Section~\ref{essn} that we can keep track of the procedures used.
We shall say that one simple extension $F$ has been obtained from
another $E$ \emph{by a standard construction} if it is the
simplification of a graph obtained by performing a sequence of
boundary outsplittings to $E$. The graph $F$ in
Theorem~\ref{essclass} has been obtained from both $SE_1$
and
$SE_2$ by a standard construction.
\end{remark}
The next example shows that the hypothesis of essentiality in
Theorem~\ref{essclass} cannot be completely dropped.
\begin{example}
Consider the following graph $G$ 
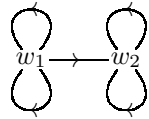
\begin{equation*}
\objectmargin{1pt}
  \xygraph{!~:{@{-}|@{>}}
{w_{1}} (:@(ur,ul)"{w_{1}}", :[r]{w_{2}} (:@(ur,ul)"{w_{2}}",
     :@(dr,dl)"{w_{2}}") , :@(dr,dl)"{w_{1}}")}
\end{equation*}
and its extensions $E_1$ and
$E_2$; 
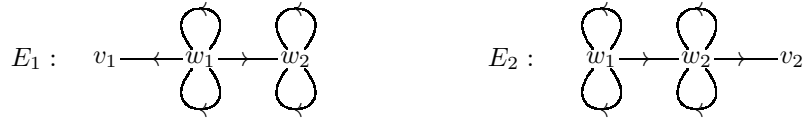
\begin{equation*}
\objectmargin{1pt}
  E_{1}:\quad\xygraph{!~:{@{-}|@{>}}
{w_{1}} (:@(ur,ul)"{w_{1}}", :[r]{w_{2}} (:@(ur,ul)"{w_{2}}",
     :@(dr,dl)"{w_{2}}") , :@(dr,dl)"{w_{1}}") :[l]{v_{1}}} 
\qquad\qquad\qquad E_{2}:\quad\xygraph{!~:{@{-}|@{>}}
{w_{1}} (:@(ur,ul)"{w_{1}}", :[r]{w_{2}} (:@(ur,ul)"{w_{2}}",
     :@(dr,dl)"{w_{2}}",:[r]{v_{2}}) , :@(dr,dl)"{w_{1}}")}
\end{equation*}
Note that $E_2$ is essential but $E_1$ is not.
On one hand, we have $A_G = \left( \begin{smallmatrix} 2 &
1 \\ 0 & 2 \end{smallmatrix} \right)$, $W_{E_1} = \left(
\begin{smallmatrix} 1 \\ 0 \end{smallmatrix} \right)$,
and $W_{E_2} = \left( \begin{smallmatrix} 0 \\ 1
\end{smallmatrix} \right)$, so
\[
W_{E_1} - W_{E_2} = \begin{pmatrix} 1 \\ -1
\end{pmatrix} = \begin{pmatrix} 1 & 1 \\ 0 & 1 \end{pmatrix}
\begin{pmatrix} 2 \\ -1 \end{pmatrix} = (A_G-I) \begin{pmatrix}
2 \\ -1 \end{pmatrix}.
\]
On the other hand, we claim that $C^*(E_1)$ is not Morita
equivalent to $C^*(E_2)$, so that they cannot have a common full
corner. To see this, recall from \cite[Theorem~4.4]{BPRS} that
the ideals in $C^*(E_i)$ are in one-to-one correspondence with
the saturated hereditary subsets of $E_i^0$. The saturated
hereditary subsets of
$E_1^0$ are
$\{ v_1 \}$, $\{ v_1, w_2 \}$, $\{ v_1, w_1, w_2
\}$ and $\{ w_2 \}$, and those of
$E_2^0$ are $\{ v_2 \}$, $\{ v_2, w_2 \}$ and $\{ v_2, w_1, w_2
\}$. Thus $C^*(E_1)$ has more ideals than $C^*(E_2)$. But if
they were Morita equivalent, the Rieffel correspondence would
set up a bijection between their ideal spaces.
\end{example}
This example shows that the way the sinks $v_i$ are attached
to $G$ can affect how the ideal $I(v_0)$ lies in the ideal space
of
$C^*(E)$. In the next section, we give a simple condition on the
way
$v_i$ are attached which ensures that the primitive ideal spaces
of $C^*(E_i)$ are homeomorphic, and show that under this
condition there is a good analogue of Theorem~\ref{essclass}.
However, there is one situation in which essentiality is not
needed: when $C^*(G)$ is an $AF$-algebra.
\begin{corollary}
Let $G$ be a graph with no sources for which $C^*(G)$ is
an $AF$-algebra, and let $(E_1,v_1)$ and
$(E_2,v_2)$ be 1-sink extensions of
$G$.  If there exists $n \in \bigoplus_{G^0} \Z$ such that
$W_{E_1} = W_{E_2} + (A_G-I)n$, then then there is a simple
1-sink extension $F$ of $G$ such that $C^*(F)$ is
$C^*(G)$-embeddable in both $C^*(E_1)$ and $C^*(E_2)$.
\label{AF}
\end{corollary}
\begin{proof}
We first recall from \cite[Theorem~2.4]{KPR} that $C^*(G)$ is
$AF$ if and only if $G$ has no loops. Now we proceed as in the
proof of Theorem~\ref{essclass}. Everything goes the same until
we come to consider the finite subgraph $D$ associated to the
support of the vector $n$. Since there are no loops in $G$, $D$
must have a sink, and the argument in the second paragraph of
the proof of Theorem~\ref{essclass} suffices; this does not use
essentiality.
\end{proof}
\section{A classification  for non-essential 1-sink
extensions}\label{iness1}
 Recall from
\cite[\S6]{BPRS} that a
\emph{maximal tail} in a graph
$E$ is a nonempty subset of $E^0$ which is cofinal under $\geq$,
is backwards hereditary ($v\geq w$ and $w\in \gamma$ imply $v\in
\gamma$), and contains no sinks (for each $w\in \gamma$, there
exists $e\in E^1$ with $s(e)=w$ and $r(e)\in \gamma$). 
\begin{definition}
Let $(E,v_0)$ be a 1-sink extension of $G$. The \emph{closure}
of the sink $v_0$ is the set
\[
\overline{v_0}:=\bigcup\{\gamma:\gamma \text{ is a maximal
tail in $G$ and }\gamma\geq v_0\}.
\]
\end{definition}
To motivate this definition, we notice first that the
extension is essential if and only if
$\overline{v_0}=G^0$. More generally (although it is not
logically necessary for our results), we
explain how this notion of closure is related to  the closure
of sets in
$\Prim C^*(E)$, as described in
\cite[\S6]{BPRS}. For each
sink $v$, let $\lambda_v:=\{w\in E^0:w\geq v\}$, and let
\[
\Lambda_E:=\{\text{maximal tails in $E$}\}\cup\{\lambda_v:v
\text{ is a sink in $E$}\}.
\]
The set $\Lambda_E$ has a topology in which the closure of a subset $S$
is $\{\lambda:\lambda\geq\bigcup_{\chi\in S}\chi\}$, and it is
proved in \cite[Corollary 6.5]{BPRS} that when $E$ satisfies
Condition~(K) of \cite{KPRR}, $\lambda\mapsto
I(E^0\setminus\lambda)$ is a homeomorphism of $\Lambda_E$ onto
$\Prim C^*(E)$. If $(E,v_0)$ is a 1-sink extension of $G$, then
the only loops in $E$ are those in $G$, so $E$ satisfies
Condition~(K) whenever
$G$ does. A subset of $G^0$ is a maximal tail in $E$ if and only
if it is a maximal tail in $G$, and because every sink in $G$ is
a sink in $E$, we deduce that
$\Lambda_E=\Lambda_G\cup \{\lambda_{v_0}\}$.
\begin{lemma}
Suppose $G$ satisfies Condition~\textnormal{(K)}, and
$(E_1,v_1)$,
$(E_2,v_2)$ are 1-sink extensions of $G$. Then
$\overline{v_1}=\overline{v_2}$ if and only if there is a
homeomorphism $h$ of $\Prim C^*(E_1)$ onto $\Prim C^*(E_2)$ such
that
\begin{equation}\label{canonhomeo}
h(I(E_1^0\setminus\lambda))=I(E_2^0\setminus\lambda)\text{ for
$\lambda\in \Lambda_G$, and } 
h(I(E_1^0\setminus\lambda_{v_1}))=I(E_2^0\setminus\lambda_{v_2}).
\end{equation}
\end{lemma}
\begin{proof}
For any 1-sink extension $(E,v_0)$, the map $J\mapsto
\pi_E^{-1}(J)$ is a homeomorphism of $\Prim C^*(G)$ onto the
closed subset $\{J\in \Prim C^*(E):J\supset I(v_0)\}$. If
$\lambda\in \Lambda_G\subset \Lambda_E$, then
$\pi_E^{-1}(I(G^0\setminus\lambda))=I(E^0\setminus \lambda)$,
and hence $h$ is always a homeomorphism of the closed set
$\{I(E_1^0\setminus \lambda):\lambda\in \Lambda_G\}$ in
$\Prim C^*(E_1)$ onto the corresponding subset of $\Prim
C^*(E_2)$. So the only issue is whether the closures of the sets
$I(E_1^0\setminus\lambda_{v_1})$ and
$I(E_2^0\setminus\lambda_{v_2})$ match up. But 
\[
\overline{I(E_i^0\setminus\lambda_{v_i})}=
\{I(E_i^0\setminus\lambda):\lambda\geq \lambda_{v_i}\}=
\{I(E_i^0\setminus\lambda):\lambda\geq v_i\}.
\]
Since other sets $\lambda_v$ associated to sinks are never
$\geq v_i$, the ideals on the right-hand side are those
associated to the maximal tails lying in
$\overline{v_i}$, and the result follows. 
\end{proof}
We now return to the problem of proving analogues of
Theorem~\ref{essclass} for non-essential extensions.  Notice
that the closure is a subset of
$G^0$ rather than
$E^0$: we have defined it this way because we want to compare
the closures in different extensions.
\begin{proposition}\label{closureclass}
Suppose that $(E_1,v_1)$ and $(E_2,v_2)$ are 1-sink extensions of
$G$ with finitely many boundary vertices, and suppose that
$\overline{v_1}=\overline{v_2}=C$, say. If
$W_{E_1}-W_{E_2}$ has the form $(A_G-I)n$ for some $n\in
\bigoplus_C\Z$, then there is a simple 1-sink extension $F$ of
$G$ such that $C^*(F)$ is $C^*(G)$-embeddable in both $C^*(E_1)$
and $C^*(E_2)$.
\end{proposition}
We aim to follow the proof of Theorem~\ref{essclass}, so we need
to check that the operations used there will not affect the
hypotheses in Proposition~\ref{closureclass}.  It is obvious
that the closure is unaffected by simplifications. It is true
but not so obvious that it is unaffected by boundary
outsplittings:
\begin{lemma}\label{outsplclosure}
Suppose $(E,v_0)$ is a 1-sink extension of a graph $G$, and $e$
is a boundary edge in $E$. Then
the closures of
$v_0$ in
$E$ and
$E(e)$ are the same.
\end{lemma}
\begin{proof}
Suppose $\gamma$ is a maximal tail such that $\gamma\geq v_0$ in
$E(e)$ and $z\in \gamma$; we want to prove $z\geq B_E^0$. We
know $z\geq w$ for some $w\in B_{E(e)}^0$. If $w\in B_E^0$,
there is no problem. If $w\notin B_E^0$, then $w=s(f)$ for some
$f\in G^1$ with $r(f)=s(e)$, so $z\geq w\geq s(e)\in B_E^0$.
Now suppose $\gamma\geq v_0$ in $E$ and $z\in \gamma$; we want
to prove that $z\geq B_{E(e)}^0$. We know that there is a path
$\alpha$ with $s(\alpha)=z$ and $r(\alpha)\in B_E^0$. If
$r(\alpha)\not= s(e)$, we have $z\geq r(\alpha)\in B_{E(e)}^0$.
If $r(\alpha)=s(e)$ and $|\alpha|\geq 1$, we have $z\geq
r(\alpha_{|\alpha|-1})\in B_{E(e)}^0$. The one remaining
possibility is that $z=s(e)$ and there is no path of length at
least 1 from $s(e)$ to $s(e)$. Because $\gamma$ is a tail, there
exists $f\in G^1$ such that $s(f)=s(e)$ and $r(f)\in\gamma$. Now
we use $\gamma\geq v_0$ to get a path $\beta$ with
$s(\beta)=r(f)$ and $r(\beta)\in B_{E^0}^0\setminus \{s(e)\}$,
and we are back in the first case with $\alpha=f\beta$.
\end{proof}
\begin{proof}[Proof of Proposition~\ref{closureclass}]
Since the closures $\overline{v_1}$ and
$\overline{v_2}$ are unaffected by simplification and boundary
outsplitting, we can run the argument of Theorem~\ref{essclass}.
In doing so, we never have to leave the common closure $C$: by
hypothesis,
$n=\sum_{k=1}^m a_k\delta_{w_k}$ for some $w_k\in C$, so all the vertices
on the subgraph $D$ used in the inductive step lie in $C$. When
$D$ has a sink, the argument goes over verbatim. When $D$ has a
loop $\gamma$, all the vertices on $\gamma$ lie in $C$, and the
hypothesis $\overline{v_1}=C=\overline{v_2}$ implies that
$\gamma\geq B_{E_1}^0\geq B_{E_2}^0\geq B_{E_1}^0$, so we can
still apply Lemma~\ref{Iainsideaforalemma}. The rest of the
argument carries over.
\end{proof}
The catch in Proposition~\ref{closureclass} is that the
vector $n$ is required to have support in the common
closure $C$. For our applications to $n$-sink extensions in
the next section, this is just what we need. However, if we
are only interested in 1-sink extensions, this requirement
might seem a little unnatural. So it is interesting that
we can often remove it:
\begin{lemma}\label{suppinC}
Suppose that $(E_1,v_1)$ and $(E_2,v_2)$ are 1-sink extensions of
$G$, and suppose that
$\overline{v_1}=\overline{v_2}=C$, say. Suppose that 1 is not an
eigenvalue of the $(G^0\setminus C)\times (G^0\setminus C)$
corner of $A_G$.
Then if $W_{E_1}-W_{E_2}$ has the form $(A_G-I)n$
for some
$n\in
\bigoplus_{G^0}\Z$, we have
$n\in\bigoplus_{C}\Z$. 
\end{lemma}
\begin{proof}
Since the maximal tails comprising $C$ are backwards hereditary,
there are no paths from $G^0\setminus C$ to $C$. Thus $A_G$
decomposes with respect to the decomposition $G^0=(G^0\setminus
C)\cup C$ as
$A_G =
\left(
\begin{smallmatrix} B&0\\ C&D
\end{smallmatrix}
\right) $, and $A_G-I = \left( \begin{smallmatrix} B-I&0\\ C&D-I
\end{smallmatrix} \right) $. Writing $n$ as $\left(
\begin{smallmatrix} k\\ m \end{smallmatrix} \right)$ and noting
that $W_{E_1}-W_{E_2}$ has support in $C$
shows that $(B-I)k=0$, which by the hypothesis on $A_G$ implies
$k=0$. But this says exactly what we want. 
\end{proof}
\section{A classification for $n$-sink extensions}\label{essn}
We say that an $n$-sink extension is \emph{essential} if
$G^0\geq v_i$ for $1\leq i\leq n$.
\begin{theorem}\label{essnclass}
let $(E,v_i)$ and $(F,w_i)$ be essential $n$-sink extensions of
$G$ with finitely many boundary vertices. Suppose that the
Wojciech vectors satisfy
\begin{equation}\label{Wvecscobound}
W_{(E;v_i)}-W_{(F;w_i)}\in
(A_G-I)\big(\textstyle{\bigoplus_{G^0}}\Z\big) \text{ for $1\leq
i\leq n$}.
\end{equation}
Then there is a simple $n$-sink extension $D$ of $G$ such that
$C^*(D)$ is $C^*(G)$-em\-bed\-da\-ble in both $C^*(E)$ and $C^*(F)$.
\end{theorem}
We shall prove this theorem by induction on $n$. At a key point
we need to convert $(n-1)$-sink
extensions to $n$-sink extensions. If $m\in
\prod_{G^0}\N$ and $(E,v_i)$ is an $(n-1)$-sink extension, we
denote by
$(E*m,v_i)$ the $n$-sink extension of $G$  obtained by adding
an extra vertex $v_n$ and $m(w)$ edges from each vertex $w\in
G^0$ to $v_n$. Note that $E*m$ has one new Wojciech vector
$W_{(E*m;v_n)}=m$, and the other Wojciech vectors are unchanged.
If $E$ is a simple extension, then so is $E*m$. Conversely, if
$(F,w_i)$ is a simple $n$-sink extension, then $F\setminus
w_n:=(F^0\setminus\{w_n\},F^1\setminus r^{-1}(w_n))$ is a simple
$(n-1)$-sink extension for which $(F\setminus w_n)*W_{(F;w_n)}$
can be naturally identified with $F$. 
We need to know how the operation $E\mapsto E*m$ interacts with our other
constructions:
\begin{lemma}\label{adjoinW}
 If $e$ is a boundary edge for $E$, then $e$ is a boundary
edge for $E*m$, and the boundary outsplittings satisfy $E(e)*m=(E*m)(e)$.
The simplification construction $E\mapsto
SE$ satisfies $S(E*m)=(SE)*m$.
\end{lemma} 
\begin{proof}
The only edges which are affected in forming $E(e)$ are $e$ and
the edges $f$ with $r(f)=s(e)$. Since none of the new edges in
$E*m$ have range in $E$, they are not affected by the
outsplitting. Simplifying collapses paths which end at one of
the sinks $v_i$, and forming $E*m$ adds only paths of length 1
ending at $v_n$, so there is nothing extra to collapse in
simplifying
$E*m$.  
\end{proof}
\begin{proof}[Proof of Theorem~\ref{essnclass}]
As in the 1-sink case, it suffices by
Proposition~\ref{simpembeds} to prove the result when $E$ and
$F$ are simple. So we assume this. Our proof is by induction on
$n$, but we have to be careful to get the right inductive
hypothesis. So we shall prove that by performing $n$ standard
constructions on both $E$ and $F$, we can arrive at simple
$n$-sink extensions of $G$ with all their Wojciech vectors equal;
these graphs are then isomorphic, and we can take $D$ to be
either of them.  Theorem~\ref{essclass} says that this is true
for $n=1$ (see Remark~\ref{standcons}). 
So we suppose that our inductive hypothesis holds for all
simple $(n-1)$-sink extensions satisfying the hypotheses of
Theorem~\ref{essnclass}. Then
$E\setminus v_n$ and
$F\setminus w_n$ are simple $(n-1)$-sink extensions of $G$ with
Wojciech vectors
$W_{(E\setminus v_n;v_i)}=W_{(E;v_i)}$ and $W_{(F\setminus
w_n;w_i)}=W_{(F;w_i)}$ for $i\leq n-1$. So the Wojciech vectors
of $E\setminus v_n$ and $F\setminus w_n$ satisfy the hypothesis
(\ref{Wvecscobound}). Since $G^0\geq v_i$ in $E$, and we
have not deleted any edges except those ending at $v_n$ and
$w_n$, we still have $G^0\geq v_i$ in $E\setminus v_n$ for $i\leq
n-1$, and similarly $G^0\geq w_i$ in $F\setminus w_n$. By the
inductive hypothesis, therefore, we can perform $(n-1)$ standard
constructions on each of $E$ and $F$ to arrive at
the same simple $(n-1)$-sink extension $(D,u_i)$ of $G$. 
By Lemma~\ref{adjoinW}, $D*W_{(E;v_n)}$ and
$D*W_{(F;w_n)}$ are obtained from $E=(E\setminus
v_n)*W_{(E;v_n)}$ and $F=(F\setminus
w_n)*W_{(F;w_n)}$ by $(n-1)$ 
standard constructions. We now view $(D^E,v_n):=D*W_{(E;v_n)}$ and
$(D^F,w_n):=D*W_{(F;w_n)}$ as two simple 1-sink extensions of the
graph $D$. Since the standard constructions have not affected
the path structure of
$G$ inside $D$, and we assumed $G^0\geq v_n$ in $E$, we still
have $G^0\geq v_n$ in $D^E$, and similarly $G^0\geq w_n$ in
$D^F$. Because any sink in $G$ has to be a sink in $E$, the
hypothesis $G^0\geq v_n$ in $E$ implies that $G$ has no sinks;
thus every vertex in $G$ lies on an infinite path $x$, and
hence in the maximal tail $\gamma:=\{v:v\geq x\}$. Thus
$G^0\geq v_n$ says precisely that $G^0$ is the closure of
$v_n$ in $D^E$. Of course the same is true of $w_n$ in $D^F$.
We can therefore apply Proposition~\ref{closureclass} to deduce
that we can by one more standard construction on each of $D^E$
and $D^F$ reach the same 1-sink extension $(C,u_n)$ of $D$;
since all the boundary vertices of $D$ lie in $G$, this standard
construction for extensions of $D$ is a also standard for
extensions of $G$, and hence $C$ can also be obtained by
performing $n$ standard constructions to each of $E$ and $F$. 
This completes the proof of the inductive hypothesis, and hence
of the theorem.
\end{proof}
\section{$K$-theory of 1-sink extensions}
\begin{proposition}
Suppose that $G$ is a row-finite graph with no sinks, and
$(E,v_0)$ is a 1-sink extension of $G$ such that $W_E\perp
\ker(A_G^t-I)$. If $(F,w_0)$ is a 1-sink extension of $G$ and
$\phi:C^*(F)\to C^*(E)$ is a $C^*(G)$-embedding, then 
there exists $k\in \prod_{G^0}\Z$ such that $W_E-W_F=(A_G-I)k$.
\label{converse}
\end{proposition}
For the proof, we need to know the $K$-theory of the
$C^*$-algebras of graphs with sinks, which was was calculated in
\cite[\S3]{RS}. We summarise some results from \cite {RS} in a
convenient form:
\begin{lemma}\label{RSKthy}
Suppose $G$ has no sinks and $(E,v_0)$ is a 1-sink extension of
$G$ with graph algebra $C^*(E)=C^*(s_e,p_v)$. Let $\psi^E$ be the
homomorphism of 
$\big(\bigoplus_{G^0}\Z\big)\oplus\Z$ into $K_0(C^*(E))$ which
is determined on the standard basis elements by
$\psi^E(\delta_v,0):=[p_v]$ for $v\in G^0$ and
$\psi^E(0,1)=[p_{v_0}]$. Then $\psi^E$ induces an isomorphism of
the cokernel of $((A_G^t-I)\oplus
W_E^t):\bigoplus_{G^0}\Z\to\big(\bigoplus_{G^0}\Z\big)\oplus\Z$
onto $K_0(C^*(E))$.
\end{lemma}
\begin{proof}
We first suppose that $(E,v_0)$ is simple. Then
$\big(\bigoplus_{G^0}\Z\big)\oplus\Z$ is the group
$\Z^{G^0}\oplus\Z^W$ considered in \cite[\S3]{RS}, and it
suffices to show that $\psi^E$ is the homomorphism
$\overline{\phi}$ considered there. To do this, we need to check
that the map $S$ of $K_0(C^*(E)\times_\gamma \T)$ onto
$K_0(C^*(E))$ in \cite[(3.3)]{RS} satisfies
$S([p_v\chi_1])=[p_v]$. The map $S$ is built up from the
homomorphisms induced by  the embedding of
$C^*(E)\times_\gamma\T$ in the dual crossed product
$(C^*(E)\times_\gamma\T)\times_{\widehat\gamma}\Z$, the
Takesaki-Takai duality isomorphism
$(C^*(E)\times_\gamma\T)\times_{\widehat\gamma}\Z\cong
C^*(E)\otimes\K(\ell^2(\Z))$, and the map $a\mapsto a\otimes p$
of
$C^*(E)$ into $C^*(E)\otimes\K(\ell^2(\Z))$ determined by a
rank-one projection $p$. The
formulas at the start of the proof of \cite[Theorem~6]{R}
show that, because
$p_v$ is fixed under
$\gamma$, the duality isomorphism carries $p_v\chi_1\in
C^*(E)\times_\gamma\T\subset
(C^*(E)\times_\gamma\T)\times_{\widehat\gamma}\Z$ into
$p_v\otimes M(\chi_1)$, where $M(\chi_1)$ is the projection onto
the subspace spanned by the basis element $e_1$. Thus $S$ has
the required property, and the result for simple extensions now
follows from
\cite[Theorem~3.2]{RS}.
If $(E,v_0)$ is an arbitrary 1-sink extension, we consider its
simplification $SE$ and the embedding $\phi^{SE}$ of $C^*(SE)$ in
$C^*(E)$ provided by Proposition~\ref{simpembeds}, which by
\cite[Proposition~1.2]{P} induces an isomorphism $\phi^{SE}_*$
in $K$-theory. But now it is easy to check that
$\phi_*^{SE}\circ\psi^{SE}=\psi^E$, and the result follows. 
\end{proof}
We now begin the proof of Proposition~\ref{converse}. Since the
image of $\phi$ is a full corner in $C^*(E)$, it induces an
isomorphism $\phi_*$ of $K_0(C^*(F))$ onto $K_0(C^*(E))$ (by,
for example, \cite[Proposition~1.2]{P}). The properties of the
$C^*(G)$-embedding $\phi$ imply that
$\phi_*([p_{w_0}])=[p_{v_0}]$ and $(\pi_E)_*\circ
\phi_*=(\pi_F)_*$. We need to know how $\phi_*$ interacts
with the descriptions of $K$-theory provided by
Lemma~\ref{RSKthy}.
\begin{lemma}\label{1stvbleok}
The induced homomorphism $\phi_*:K_0(C^*(F))\to K_0(C^*(E))$
satisfies $\phi_*(\psi^F(0,1))=\psi^E(0,1)$, and for each $z\in
\bigoplus_{G^0}\Z$, there exists $\ell\in \Z$ such that
$\phi_*(\psi^F(z,0))=\psi^E(z,\ell)$.
\end{lemma}
\begin{proof}
The first equation is a translation of the condition
$\phi_*([p_{w_0}])=[p_{v_0}]$. For the second, let
$\psi^G:\bigoplus_{G^0}\Z\to K_0(C^*(G))$ be the 
homomorphism such that $\psi^G(\delta_v)=[p_v]$, which
induces the usual isomorphism of $\coker (A_G^t-I)$ onto
$K_0(C^*(G))$. If $\rho:(\bigoplus_{G^0}\Z)\oplus\Z\to \bigoplus_{G^0} \Z$ is
given by
$\rho(z,\ell):=z$, then we have
$(\pi_E)_*\circ\psi^E=\psi^G\circ\rho$, and similarly for $F$.
Thus
\begin{equation}\label{effectphi*}
(\pi_E)_*\circ\phi_*\circ\psi^F=(\pi_F)_*\circ\psi^F
=\psi^G\circ\rho.
\end{equation}
Now fix $z\in \bigoplus_{G^0}\Z$. Since $\psi^E$ is surjective,
there exists $(x,y)\in \big(\bigoplus_{G^0}\Z\big)\oplus\Z$ such
that $\psi^E(x,y)=\phi_*(\psi^F(z,0))$. From (\ref{effectphi*})
we have
\[
\psi^G(z)=(\pi_E)_*\circ\phi_*\circ\psi^F(z,0)
=(\pi_E)_*\circ\psi^E(x,y)=\psi^G(x),
\]
and hence there exists $u\in\bigoplus_{G^0}\Z$ such that
$x=z+(A_G^t-I)u$. Now because $\psi^E$ is constant on the image
of $(A_G^t-I)\oplus
W_E^t$, we have
\[
\phi_*(\psi^F(z,0))=\psi^E(x,y)=\psi^E(z+(A_G^t-I)u,y)=
\psi^E(z,y-W_E^tu),
\]
and $\ell:=y-W_E^tu$ will do.
\end{proof}
\begin{proof}[Proof of Proposition~\ref{converse}]
By Lemma~\ref{1stvbleok}, for each $v\in G^0$ there exists
$k_v\in\Z$ such that
$\phi_*(\psi^F(\delta_v,0))=\psi^E(\delta_v,k_v)$. We define
$k=(k_v)\in
\prod_{G^0}\Z$. A calculation shows that for any $(y,\ell)\in
\big(\bigoplus_{G^0}\Z\big)\oplus\Z$ we have
\begin{align}\label{genform}
\phi_*\circ\psi^F(y,\ell)&=\sum_vy_v(\phi_*\circ\psi^F)
(\delta_v,0)+\ell(\phi_*\circ\psi^F)(0,1)\\
&\notag=\Big(\sum_v\psi^E(y_v\delta_v,y_vk_v)\Big)+\ell\psi^E(0,1)\\
&\notag=\psi^E(y,k^ty+\ell).
\end{align}
Now let $z\in \bigoplus_{G^0}\Z$. On one hand, we have from
(\ref{genform}) that
\begin{equation}\label{onehand}
\phi_*\circ\psi^F\big(((A_G^t-I)\oplus W_F^t)(z)\big)
=\psi^E\big((A_G^t-I)z,k^t(A_G^t-I)z+W_F^tz\big).
\end{equation}
On the other hand, since $\psi^F\circ((A_G^t-I)\oplus W_F^t)=0$,
its composition with $\phi_*$ is also $0$. Thus the class (\ref{onehand})
must vanish in $K_0(C^*(E))$, and there exists $x\in \bigoplus_{G^0}\Z$
such that 
\begin{equation}\label{otherhand}
\big((A_G^t-I)z,k^t(A_G^t-I)z+W_F^tz\big)
=\big((A_G^t-I)x,W_E^tx\big).
\end{equation}
Comparing (\ref{onehand}) and (\ref{otherhand}) shows that
$x-z\in \ker (A_G^t-I)$ and 
\[
k^t(A_G^t-I)z+W_F^tz=W_E^tx=W_E^tz+W_E^t(x-z).
\]
Since we are supposing $W_E\perp\ker(A_G^t-I)$, we deduce that
$W_E^t(x-z)=0$. We have now proved that 
\[
k^t(A_G^t-I)z=(W_E^t-W_F^t)z\ \mbox{ for all $z\in
\textstyle{\bigoplus_{G^0}}\Z$},
\]
 which implies $(A_G-I)k=W_E-W_F$, as required.
\end{proof}
\begin{corollary}
Suppose that $G$ is a row-finite graph with no sinks and with
the property that $\ker (A_G^t-I) = \{ 0 \}$.  Let
$(E_1,v_1)$ and $(E_2,v_2)$ be 1-sink extensions of $G$. 
If there is a 1-sink extension $F$ such that $C^*(F)$ is
$C^*(G)$-embeddable in both $C^*(E_1)$ and $C^*(E_2)$,  then
there exists
$k\in
\prod_{G^0}\Z$ such that 
$W_{E_1}-W_{E_2} =(A_G-I)k$.
\label{corconverse}
\end{corollary}
\begin{corollary}
Suppose that $G$ is a finite graph with no sinks or sources
whose vertex matrix $A_G$ satisfies $\ker(A_G^t-I)=0$. Let
$(E_1,v_1)$ and $(E_2,v_2)$ be 1-sink extensions of $G$ such
that $\overline{v_1}=\overline{v_2}$. Then there is a 1-sink
extension $F$ such that $C^*(F)$ is
$C^*(G)$-embeddable in both $C^*(E_1)$ and $C^*(E_2)$ if and
only if there exists
$k\in
\bigoplus_{G^0}\Z$ such that 
$W_{E_1}-W_{E_2} =(A_G-I)k$.
\end{corollary}
\begin{proof}
The forward direction follows from the previous corollary. For
the converse, we seek to apply Proposition~\ref{closureclass}.
To see that $n$ has support in the common closure
$C:=\overline{v_1}=\overline{v_2}$, recall that $A_G$
decomposes  as
$A_G =
\left(
\begin{smallmatrix} B&0\\ C&D
\end{smallmatrix}
\right) $ with respect to $G^0=(G^0\setminus
C)\cup C$. Thus $1$ is an eigenvalue for the $(G^0\setminus
C)\times (G^0\setminus
C)$ corner $B$ of $A_G$ if and only if it is an eigenvalue for
$A_G$, and hence if and only if it is an eigenvalue for
$A_G^t$. So Lemma~\ref{suppinC} applies, $\supp n$ lies in $C$,
and the result follows from Proposition~\ref{closureclass}.
\end{proof}
\bibliographystyle{ams}
\providecommand{\bysame}{\leavevmode\hbox to3em{\hrulefill}\thinspace}

\end{document}